\documentclass[11pt,twoside]{article}%

\usepackage[latin1]{inputenc}%
\usepackage{amsmath,amssymb,latexsym}%
\usepackage{theorem}
\renewcommand{\emph}{\textbf}%
\def\endrem{}%

\def\proof#1{{\def\temp{#1} 
              \ifx\temp\empty 
                  \noindent\slshape\textbf{proof:\ }
              \else 
                  \noindent\slshape\textbf{proof of\ #1:\ }
              \fi}}

\def\qed{\hspace*{\fill}$\square$}

\newcommand{\ol}{\overline}%
\newcommand{\wt}{\widetilde}%
\newcommand{\normaleq}{{\,\raisebox{.17ex}{\§$\vartriangleleft$\§}
    \mkern-19mu\leqslant}}%
\def\Gn#1_#2.{\def\temp{#1} \mathbf{G}_{#2}%
      \ifx\temp\empty\else (#1) \fi}%
\def\G#1.{\Gn#1_.}%
\def\Gan#1.{\def\temp{#1} \mathbf{G}^\mathsf{an}%
           \ifx\temp\empty\else (#1) \fi}%
\def\Gis#1.{\def\temp{#1} \mathbf{G}^\mathsf{is}%
           \ifx\temp\empty\else (#1) \fi}%
\def\AdGan#1.{\def\temp{#1} \Ad\mathbf{G}^\mathsf{an}%
           \ifx\temp\empty\else (#1) \fi}%
\def\AdGis#1.{\def\temp{#1} \Ad\mathbf{AdG}^\mathsf{is}%
           \ifx\temp\empty\else (#1) \fi}%
\def\AdG#1.{\def\temp{#1} \mathrm{Ad}\mathbf{G}%
           \ifx\temp\empty\else (#1) \fi}%
\def\wtG#1.{\def\temp{#1} \wt{\mathbf{G}}%
           \ifx\temp\empty\else (#1) \fi}%
\let\jura=\S %
\def\S#1.{\def\temp{#1} \mathbf{S}%
           \ifx\temp\empty\else (#1) \fi}%
\def\wtS#1.{\def\temp{#1} \wt{\mathbf{S}}%
           \ifx\temp\empty\else (#1) \fi}%
\def\T#1.{\def\temp{#1} \mathbf{T}%
           \ifx\temp\empty\else (#1) \fi}%
\def\M#1.{\def\temp{#1} \mathbf{M}%
           \ifx\temp\empty\else (#1) \fi}%
\def\P#1_#2.{\def\temp{#1} \mathbf{P}_{#2}%
              \ifx\temp\empty\else (#1) \fi}%
\def\olP#1_#2.{\def\temp{#1} \ol{\mathbf{P}}_{#2}%
           \ifx\temp\empty\else (#1) \fi}%
\def\Ps#1_#2^#3.{\def\temp{#1} \mathbf{P}_{#2}^{#3}%
           \ifx\temp\empty\else (#1) \fi}%
\def\U#1_#2.{\def\temp{#1} \mathbf{U}_{#2}%
           \ifx\temp\empty\else (#1) \fi}%
\def\Z#1_#2:#3.{\def\Ztemp{#3} \mathcal{Z}_{#2}(#1)%
           \ifx\Ztemp\empty\else (#3) \fi}%
\def\N#1_#2:#3.{\def\Ntemp{#3} \mathcal{N}_{#2}(#1)%
           \ifx\Ntemp\empty\else (#3) \fi}%
\def\Ru#1:#2.{\def\Rutemp{#2}  R_u(#1)%
           \ifx\Rutemp\empty\else (#2) \fi}%
\def\Trans#1#2_#3.{\mathsf{Trans}_{#3}(#1,\,#2)}%
\def\Aut#1_#2.{\mathsf{Aut}_{#2}(#1)}%
\def\Stab#1_#2.{#2_{\{#1\}}}%
\def\Fix#1_#2.{{#2_{#1}}}
\newcommand{\ZZ}{\mathbb{Z}}%
\newcommand{\Ad}{\operatorname{Ad}}%
\newcommand{\id}{{\sf id}}%
\newcommand{\modul}{\mathsf{mod}}%
\newcommand{\im}{\operatorname{im}}%
\def\Romannumeral#1 {\uppercase\expandafter{\romannumeral#1}}%
\def\§{\hspace*{-\the\mathsurround}}%
\newcommand{\too}{\mathop{-\!-\!\!\!\longrightarrow}\limits}%
\def\Splited#1>#2#3<#4{\unitlength.65mm%
\mbox{$\displaystyle #1$\rule[-3.6mm]{0mm}{3.6mm}%
\begin{picture}(25,10)(0,5.5)
      \put(2.5,7){\vector(1,0){20}}
      \put(12,8.5){$\scriptstyle #2$}
      \put(17.5,5){\oval(10,10)[br]}
      \put(17.5,0){\vector(-3,1){15}}
      \put(12.5,2){$\scriptstyle #3$}
\end{picture}
$\displaystyle #4$}}
\def\picArrow#1{\unitlength.65mm%
\begin{picture}(25,10)(0,5.5)
      \put(2.5,7){\vector(1,0){20}}
      \put(12,8.5){$\scriptstyle #1$}
\end{picture}}
\def\SES#1#2>#3>#4#5<#6{%
\mbox{$\displaystyle 1$\picArrow{}
$\displaystyle #1$\picArrow#2\Splited#3>#4#5<{#6}
\picArrow{}$\displaystyle 1$}}
\def\picarrow#1{\unitlength.35mm%
\begin{picture}(25,10)(0,9.3)
      \put(2.5,13){\vector(1,0){20}}
      \put(12,15){$\scriptstyle #1$}
\end{picture}}
\def\ses#1#2>#3>#4#5<#6{%
\mbox{$\textstyle 1$\picarrow{}
$\textstyle #1$\picarrow#2\splited#3>#4#5<{#6}
\picarrow{}$\textstyle 1$}}
\def\splited#1>#2#3<#4{\unitlength.35mm%
\mbox{$\textstyle #1$\rule[-3.00mm]{0mm}{3.00mm}%
\begin{picture}(25,10)(0,9.3)
      \put(2.5,13){\vector(1,0){20}}
      \put(12,15){$\scriptstyle #2$}
      \put(17.5,11){\oval(10,10)[br]}
      \put(17.5,6){\vector(-3,1){15}}
      \put(17.5,0){$\scriptstyle #3$}
\end{picture}
$\textstyle #4$}}
\begingroup%
\theoremstyle{plain}\theorembodyfont{\slshape}%
\newtheorem{satz}{Proposition}[section]%
\newtheorem{corollar}[satz]{Corollary}%
\newtheorem{lemma}[satz]{Lemma}%
\newtheorem{bemerkung}[satz]{{Remark}}%
\newtheorem{notation}[satz]{{Notation}}%
\newtheorem{definition}[satz]{{Definition}}%
\newtheorem{beispiel}[satz]{{Example}}
\endgroup%

\begingroup%
\theoremstyle{plain}
\theorembodyfont{\itshape}\theoremheaderfont{\bfseries}%
\newtheorem{theorem}[satz]{Theorem}%
\endgroup%

\begingroup%
\theoremstyle{break}\theorembodyfont{\slshape}%
\newtheorem{itembemerkung}[satz]{{Remark}}%
\newtheorem{itembeispiel}[satz]{{Example}}%
\newtheorem{itemlemma}[satz]{{Lemma}}%
\endgroup%

\begingroup%
\theoremstyle{break}\theorembodyfont{\sffamily}
\newtheorem{algo}[satz]{Algorithm}%
\endgroup%

\let\variable=\epsilon\let\epsilon=\varepsilon\let\varepsilon=\variable

\long\def\weg#1{}%

\let\setminus=\smallsetminus%
\let\variable=\epsilon\let\epsilon=\varepsilon\let\varepsilon=\variable

\newcommand{\tdlcG}{{totally disconnected locally compact group}}
\newcommand{\tdlc}{{totally disconnected locally compact}}

\begin{document}%
\title{Contraction groups and scales of automorphisms of
       totally disconnected locally compact groups}

\author{Udo Baumgartner and George A. Willis}


\thanks{Research supported by DFG Grant BA 1966/1-1 and ARC Grant
A69700321}


\maketitle

\abstract
\noindent 
We study contraction groups for automorphisms of \tdlcG s using 
the scale of the automorphism as a tool.
The contraction group is shown to be unbounded when the 
inverse automorphism has non-trivial scale 
and this scale is shown to be the inverse 
value of the modular function on 
the closure of the contraction group  
at the automorphism. 
The closure of the contraction group is represented 
as acting on a homogenous tree and closed contraction groups 
are characterised. 

\section{Introduction}
Interest in contraction groups
and related concepts
has been stimulated by
applications in the theory of probability measures and random walks on
groups  and in representation theory.

In representation theory, contraction groups
bring about
the 
Mautner phenomenon.
This is manifest in Wang's examination
\cite{Mautner(p-adic.Lie)}
of the phenomenon in $p$-adic Lie groups,
though more in the background in Moore's
treatment of the Lie group case \cite{Mautner(gen.unitary-rep)}.
To our knowledge
the first to define and study
contraction groups
(in a slightly more general context than ours, see \cite{kontr.Ext+*G})
was M\"uller-R\"omer, 
who did so to study a representation theoretic question as well
(the Wiener and Tauber property of a group algebra).

When studying
semistable convolution semigroups of probability measures
on locally compact groups, contraction groups arise naturally
\cite{Aut(Lie):contr/cp<+Appl(semist*semiG),contrG+semist-meas(p-adicLie)}.
For this reason
a lot of preliminary work on contraction groups was done
by Siebert, and most of the known results are either
derived or referenced in \cite{contrAut(lcG)}.  
Proposition~4.2 in  \cite{contrAut(lcG)} 
reduces the study of locally compact contraction groups to
a separate study of the connected and totally disconnected cases.
The connected case is covered by Corollary~2.4 in {\it  loc.cit.\/}.
In some totally disconnected groups the contraction groups 
are closed and these are studied in
\cite{semistable.conv-semiG+T(contrG)}, 
see also Remark~\ref{Ex(closedU)} below. 
We contribute to the full picture by treating 
general totally disconnected groups.

Just as Lie techniques 
are used to study contraction groups in the connected case, 
so the notions of scale and tidy subgroup  
(introduced in \cite{tdlcG.structure}) 
are useful in the totally disconnected case. 
The relevant properties of tidy subgroups 
are summarized in section~\ref{basic(s,V)}. 
Their connection to contraction subgroups
and related concepts 
is explored 
in section~\ref{U&P}.
The basic properties are collected in subsection~\ref{basic(U,P)}. 
This subsection culminates in our main result, 
Theorem~\ref{gen_quotient-U}, 
whose proof rests on the assumption  
that the group is metrizable. 

Therefore, in the following sections
we assume that all groups considered are
totally disconnected locally compact metric 
unless explicitly stated otherwise. 

From the reinterpretation of the scale function as the modular
function restricted to various subgroups (Proposition~\ref{s=mod()})
we infer that contraction groups
of automorphisms whose inverse has non-trivial scale are unbounded
(Corollary~\ref{U bounded iff}).
We succeed in Proposition~\ref{U0=1 iff} to characterize closed 
contraction groups.
Theorem~\ref{T.alpha} in Section~\ref{S:tree} then
represents contraction groups as groups of automorphisms
of a homogeneous tree.


We stick to the following conventions:
$0$ is a natural number.
The relations  $\subset$, $\vartriangleleft$ {\it  etc.\/} always imply
strict inclusion. 
Any automorphism of a topological group will be assumed to be a homeomorphism. 
By ``$X$ is stable under $\alpha$'' we mean $\alpha(X)=X$
whereas ``$X$ is invariant under $\alpha$'' means $\alpha(X)\subseteq X$.
The modular function $\modul_G$ of a locally compact topological group
$G$ is defined by the equation
$\mu(Mg)=\modul_G(g)\mu(M)$ where $\mu$ is a left Haar measure on $G$.
We use $e$ for the unit element of a group and $1$ for the trivial
group.
The function $|\cdot|$ stands for the absolute value on complex or
$p$-adic numbers, as the case may be.

\noindent{\small\textsc{Acknowledgements: }\small 
A preliminary 
version of this paper 
was completed 
while both authors 
were visiting 
the Tata Institute for Fundamental Research 
in Bombay 
in 1999. 
We wish 
to thank 
the Institute 
for its warm hospitality 
and  
Riddhi Shah 
for helpful 
discussions 
on 
an earlier form 
of our Main Result 
which now 
forms 
the foundation 
for the whole paper.}

\section{Basic facts about the scale function}
\label{basic(s,V)}
%
\label{U,P:1}

Tidy subgroups for an automorphism $\alpha$ of a \tdlcG\ $G$  
provide a local description of $\alpha$. 
For a subgroup $V$ of $G$ define subgroups $V_+$, $V_-$ of $V$ 
and $V_{++}$, $V_{--}$ of $G$ by: 
$$
V_{\pm} = \bigcap_{n\geq0}\alpha^{\pm n}(V)\quad\
\text{and}\quad
V_{\pm\pm}=\bigcup_{n\geq0}\alpha^{\pm n}(V_\pm)\;.
$$
The automorphism $\alpha$ magnifies $V_+$ and shrinks $V_-$, 
that is, 
$\alpha(V_+)\geqslant V_+$ and
$\alpha^{-1}(V_-)\geqslant V_-$.
We let $V_0$ denote the ``neutral'' part:
$$
V_0:=V_+\cap V_- = \bigcap_{k\in\mathbb{Z}}\alpha^k(V) =
\bigcap_{j\in \mathbb{N}}\alpha^{j}(V_-)=
\bigcap_{j\in \mathbb{N}}\alpha^{-j}(V_+)\;.
$$
It is stable under $\alpha$.

\begin{definition}
Let $G$ be a \tdlcG\ and $\alpha$ an automorphism of $G$.
A compact-open subgroup $V$ is called \emph{tidy for $\alpha$ in $G$} (or
\emph{tidy} for short if $\alpha$ is understood) if it satisfies
\begin{itemize}
\item[\textbf{(T1)}] $V=V_+V_-\ (=V_-V_+)$
\item[\textbf{(T2)}] $V_{--}$ (and $V_{++}$) are closed
\end{itemize}
The integer
$$
s_G(\alpha):=|\alpha(V_+):V_+|=|\alpha(V):V\cap \alpha(V)|
$$
is called the \emph{scale of $\alpha$}.
The function $s_G\colon G\to \mathbb{N}$ obtained by restricting
attention to inner automorphisms will be called the
\emph{scale function} of $G$.
\end{definition}\endrem

The scale of an automorphism is well defined by 
Theorem~2 in \cite{tdlcG.structure}. 
Observe that a subgroup $V$ which is tidy for $\alpha$ 
will be tidy for $\alpha^{-1}$ as well.

Tidy subgroups exist. 
They may be found by
running the following algorithm.

\begin{algo}[{ {\it cf.\/}} proof of Theorem 3.1 in
\cite{furtherP(s(tdG))}]\label{algo} 
{[0]} 
Choose a compact open subgroup $O\leqslant G$.\hfill\\
{[1]}
Let $^kO:=\bigcap_{i=0}^k \alpha^i(O)$.
We have $^kO_+=O_+$ and 
$^kO_-
=\alpha^k(O_-)$.\hfill\\
\hphantom{{[1]}}
For some $n\in \mathbb{N}$ (hence for all $n'\ge n$) the group $^nO$ satisfies
\textbf{(T1)}.\hfill\linebreak
\hphantom{{[2]}}
Put $O':={^nO}$.\\
{[2]} 
Let 
$\mathcal{L}:=\{x\in G \colon \alpha^i(x)\in O'\
\text{for almost all}\ i\in \mathbb{Z}\}$ 
and $L:=\ol{\mathcal{L}}$.\hfill\\ 
{[3]} 
Form
$O^*:=\{x\in O'\colon lxl^{-1}\in O'L\,\,\forall l\in L\}$ 
and define $O'':=O^* L$. 
\hphantom{{[2]}}
The group $O''$ is tidy and we output $O''$.
\end{algo}


Analyzing the above algorithm
({\it cf.\/} \cite[Theorem 3.1]{furtherP(s(tdG))}) one
arives at the alternative characterization of the scale
function as a minimal distortion value. 
$$ 
s_G(\alpha)=\min\{|\alpha(O):O\cap \alpha(O)|\colon
               O \ \text{a compact, open subgroup}\}
$$


%

The following simple criterion ensuring tidiness will
be used repeatedly in the sequel.

\begin{beispiel}\label{V=V-}
A compact, open subgroup $V$ satisfying $V=V_-$ is tidy.
\qed
\end{beispiel}

It follows immediately from the definition of the scale that $\alpha$
has a tidy subgroup satisfying this criterion if and only if
$s(\alpha) = 1$.

\section{Contraction groups and parabolics}\label{U&P}
Let $\alpha$ be an automorphism of a locally compact
group~$G$. 
This section 
links the subgroups tidy for $\alpha$ 
to two other subgroups having global dynamical definitions 
in terms of $\alpha$

\subsection{Definitions and basic properties}\label{basic(U,P)}
We define the parabolic group and the contraction group 
of $\alpha$ and begin to elucidate 
links between the various groups. 
Tidy subgroup methods show 
that the parabolic group is closed. 
They are also used in the proof of 
Theorem~\ref{gen_quotient-U}, 
which is fundamental for 
all of our results on contraction groups. 


\begin{bemerkung}\label{V- in P}
If $V$ is a compact, open subgroup of $G$, then we have
$$
V_-\subseteq P_\alpha:=\left\{ x\in G\colon 
\{\alpha^n(x)\colon n\in\mathbb{N}\}\
\text{is bounded} \right\}
$$
and $V_-=V\cap P_\alpha$ if $V$ is tidy for $\alpha$
(see Lemma 9 
in \cite{tdlcG.structure}).
\end{bemerkung}\endrem

The group $P_\alpha$ is closed
(Proposition 3 parts (iii) and (ii) in \cite{tdlcG.structure} show this)
and obviously contains
the group
$$
U_\alpha:=\{ x\in G\colon \alpha^n(x)\too_{n\to\infty} e\}\;,
$$
 which need not be closed (see Example
\ref{1stEx:U,P}(2)).

\begin{notation}\label{def P,U}
We call $P_\alpha$ and $U_\alpha$ respectively the \emph{parabolic}
subgroup and the \emph{contraction} group associated to $\alpha$.
We also let $M_\alpha:=P_\alpha\cap P_{\alpha^{-1}}$ and call it the
\emph{Levi factor} attached to $\alpha$.
The Levi factor is the set of all elements of the ambient group whose
$\langle\alpha\rangle$-orbit is bounded.
If the automorphism 
$\alpha$ is inner and is conjugation by 
$g$, we relax notation and write $P_g$, $U_g$ and $M_g$.
\end{notation}

The term `parabolic group' is suggested by Example~\ref{1stEx:U,P} to
follow.
The name `contraction group' however is standard, see
\cite{kontr.Ext+*G, Mautner(p-adic.Lie),
contrAut(lcG), semistable.conv-semiG+T(contrG)}.
Results about contraction groups, notably in the case where $U_\alpha$ is
closed, may be found in these papers. 
The next result is one such. 
We shall use it frequently.

\begin{lemma}[
              \cite{Mautner(p-adic.Lie)}, Proposition 2.1]\label{cp-contr}
Let $G$ be a locally compact group and let $\alpha$ 
be an  automorphism of $G$ such that $U_\alpha=G$. 
Then $\alpha$ is compactly contractive, that is, for any compact subset
$C$ of $G$ and any neighborhood $O$ of $e$ we have
$\alpha^n(C)\subseteq O$ for all $n\ge N(C,O)$.
\end{lemma}

Since $U_\alpha$ is $\alpha$-stable, this will imply that whenever
the ambient group is locally compact and $U_\alpha$ is closed, the
restriction of $\alpha$ to $U_\alpha$ is compactly contractive.

Since $\alpha$ shrinks $V_-$, one feels that there should be an even
closer connection between $V_-$ and $U_\alpha$ than the one between
$V_-$ and $P_\alpha$ displayed by Remark \ref{V- in P}.
Before establishing this, we note \textbf{some elementary properties
of the groups $U_\alpha$, $P_\alpha$ and $M_\alpha$}:
$$
U_{\alpha^n}=U_\alpha\ \ \text{and}\ \ P_{\alpha^n}=P_\alpha\quad
\text{for all}\ n\in \mathbb{N}\setminus\{0\}\;.
$$
$$
U_{\mathrm{id}}=1\ \ \text{and}\ \ P_{\mathrm{id}}=G,\quad
\beta(U_\alpha)=U_{\beta\alpha\beta^{-1}},\quad
\beta(P_\alpha)=P_{\beta\alpha\beta^{-1}}
$$
Further, as is plain from the definitions, when computing
the contraction groups and parabolics inside a
subgroup (stable under the automorphism in question)
we get the intersections of the contraction groups and parabolics
in the ambient group respectively with the subgroup.

Obviously $M_\alpha=M_{\alpha^{-1}}$ and $V_0\leqslant V\cap M_\alpha$
for every compact open subgroup with equality if $V$ is tidy thanks
to Remark \ref{V- in P}.

\begin{satz}\label{U normal P}
Let $G$ be a locally compact group and let $\alpha$ be 
an  automorphism.
Then $U_\alpha$ is normal in $P_\alpha$, hence
$U_\alpha M_\alpha\leqslant P_\alpha$.
\end{satz}
\proof{} Let $x\in P_\alpha$, $u\in U_\alpha$ be given. 
By definition of
$P_\alpha$ the set $\ol{\{\alpha^n(x)\colon n\in \mathbb{N}\}}$,
name it $K$, is compact.
Given an open neighborhood $O$ of $e$, choose an open neighborhood
$O'$ of $e$ satisfying $\bigcup_{k\in K} kO'k^{-1}\subseteq O$
(\cite[II.4.9]{GMW115}).
Then from $\alpha^n(u)\in O'$ for all $n\ge n_0$ we infer
$\alpha^n(xux^{-1})=\alpha^n(x)\alpha^n(u)\alpha^n(x)^{-1}\in
\bigcup_{k\in K} kO'k^{-1}\subseteq O$, proving that indeed
$xux^{-1}\in U_\alpha$. \qed

Theorem II.4.9 from \cite{GMW115} can also be used to prove the
first claim in the following result, which is useful in computations.

\begin{lemma}\label{red71.2}
Let $G$ be a locally compact group and let $d,\, v\in G$ be such that
$dv=vd$ and $\langle v\rangle$ is bounded.
Then $U_{dv}=U_d$ and  $P_{dv}=P_d$.
\qed
\end{lemma}
%

We note the following simple consequence.

\begin{corollar}\label{U neq 1->}
Let $G$ be a locally compact group and $g\in G$ be given.
Then either of $U_g\neq 1$ or $P_g\neq G$ implies that
$\ol{\langle g\rangle}=\langle g\rangle$
is infinite cyclic.
\qed
\end{corollar}

Next, we investigate the behavior of contraction groups and
parabolics under quotient maps.

\begin{satz}\label{quotient-U}
Let $p\colon G\to \ol{G}$ be a homomorphism of totally disconnected 
locally compact metric  groups which is an identification.
Let $\alpha$ be an  automorphism of $G$ leaving 
the kernel of $p$
stable and thus inducing an  automorphism
$\overline{\alpha}$ of $\overline{G}$.
Then $p(U_\alpha)=U_{\ol{\alpha}}$.
\end{satz}

Proposition \ref{quotient-U} is obtained as a corollary of the next 
result taking $H:= \ker p$. 
We first introduce some notation. 
Let $H$ be a subset of the topological group $G$. 
Call a sequence of elements $(x_n)_{n\in \mathbb{N}}$ in $G$ 
convergent  to $e$ modulo $H$, if for any neighborhood $W$ of $H$ there is
an integer $N_W$ such that $W$ contains all terms of the subsequence $(x_n)_{n\ge N_W}$ and write
$\lim_{n\in\mathbb{N}} x_n = e \mod H$ in this case.
Let 
$U_{\alpha/H}:=\{x\in G\colon \lim_{n\in \mathbb{N}} \alpha^{n}(x)= e \mod H\}$.

\begin{theorem}\label{gen_quotient-U}
Let $G$ be a totally disconnected locally compact metric  group, 
$\alpha$ an  automorphism of $G$ and 
$H$ an $\alpha$-stable 
closed 
subgroup of $G$. 
Then $ U_{\alpha/H}=U_\alpha H$. 
\end{theorem}

The criterion provided by the following lemma is useful in the proof. 

\begin{lemma}\label{conv.mod.H=acc.in.H}
Let $(x_n)_{n\in \mathbb{N}}$ be a sequence 
in a locally compact group $G$ 
and let $H$ be a subset of $G$. 
Then the following statements hold. 
\begin{itemize}
\item[(1)] Let $\{x_n\colon n\in\mathbb{N}\}$ be bounded. 
           If $H$ contains 
           each accumulation point of $(x_n)_{n\in \mathbb{N}}$ 
           then $(x_n)_{n\in \mathbb{N}}$ converges 
           to $e$ modulo $H$. 
\item[(2)] If $(x_n)_{n\in \mathbb{N}}$ converges 
           to $e$ modulo $H$, then each accumulation point of 
           $(x_n)_{n\in \mathbb{N}}$ is contained in $\overline{H}$. 
\item[(3)] If $\{x_n\colon n\in\mathbb{N}\}$ is bounded and $H$ is closed 
           then $(x_n)_{n\in \mathbb{N}}$ converges to $e$ modulo $H$ 
           iff 
           $H$ contains each accumulation point of 
           $(x_n)_{n\in \mathbb{N}}$. 
\end{itemize}
\end{lemma} 

\proof{lemma} 
(2) follows from the definitions, while 
(3) is implied by the statements (1) and (2). 
It remains to prove (1). 

We argue by contradiction. 
Let $W$ be a neighborhood of $H$ such that 
its complement contains infinitely many elements of 
the set 
$\{x_n\colon n\in\mathbb{N}\}$. 
Since $\{x_n\colon n\in\mathbb{N}\}$ is bounded, 
the subsequence $(x_{n_i})_{i\in \mathbb{N}}$ of $(x_n)_{n\in \mathbb{N}}$ 
formed by the elements not in $W$ has an accumulation point, 
which is in $H$ by assumption. 
Then $W$ must contain infinitely many elements of the set 
$\{x_{n_i}\colon i\in\mathbb{N}\}$ contradicting its definition. 
\qed

\proof{theorem}
The right hand side is contained in the left hand side.
To show the opposite inclusion, we need to show that whenever
$\alpha^{n}(x)$ converges to $e$ modulo $H$, then there is an element $h$ 
of $H$, such that $\alpha^{n}(xh)$ converges to $e$.

Since $G$ is metric, there is a decreasing sequence 
$(O^{(i)})_{i= 1}^\infty$ of compact open subgroups of $G$ 
with trivial intersection. 
Put $O^{(0)}:=G$. 
We use induction on $i\in \mathbb{N}$ to show that there are sequences
of elements $(y_i)_{i\in \mathbb{N}}$ in $G$ 
and natural numbers $(N_i)_{i\in \mathbb{N}}$
such that 
\begin{itemize}
\item[(1)] $y_0=x,\quad y_{i+1}\in y_i(O^{(i)}_0\cap H)$ 
           for all $i\in\mathbb{N}$,  
\item[(2)] $\alpha^{n}(y_i)\in O^{(i)}$ for all $n\ge N_i$, 
\item[(3)] $\lim_{n\in\mathbb{N}}\alpha^{n}(y_i)= e 
            \mod O^{(i)}_0\cap H \;.$
\end{itemize} 

Putting $y_0:=x$ and $N_0:=0$ provides a basis for the induction.  
For the induction step we will use the following lemma.

\begin{lemma}
Let $G$ be a totally disconnected locally compact group, 
$\alpha$ an  automorphism of $G$, 
$H$ an $\alpha$-stable closed subgroup of $G$ 
and $x$ an element of $U_{\alpha/H}$.
Then, 
for any compact open subgroup $O$ of $G$ there is an element $h$ in $H$ and
a natural number $N$ such that $\alpha^n(xh)$ is contained in $O$ for each
$n\ge N$.
The sequence $\left(\alpha^n(xh)\right)_{n\in \mathbb{N}}$ converges to $e$ modulo 
$O_0 \cap H$. 
\end{lemma}
{\slshape\textbf{proof of lemma:\ }}
Applying step~1 of the Algorithm \ref{algo} to $O\cap H$, 
we may assume that the
intersection of $O$ with $H$ satisfies property (T1) with respect to $\alpha$,
hence $O\cap H= (O\cap H)_- (O\cap H)_+$.
Using continuity of $\alpha$, choose a compact open subgroup $V$ of $O$ such
that $\alpha(V)\subseteq O$.
Then $\alpha(V(O\cap H))\subseteq O\alpha(O\cap H)$.
Let $N$ be such that $\alpha^n(x)$ is contained in $VH$ for all $n\ge N$.

Choose an element $h_0$ in $H$ such that $\alpha^N(xh_0)\in V$.
We will complete $h_0$ recursively to a sequence $(h_i)$ of elements in
$h_0(O\cap H)$ such that
$$
\alpha^{N+j}(xh_i)\in V(O\cap H) \ \text{for}\ 0\le j\le i \;.
$$
The recursion starts with $h_0$.
Suppose then, that $h_0$ up to $h_k$ have already been constructed to
satisfy this condition.\\
Then, using $\alpha(O\cap H)\subseteq (O\cap H) \alpha((O\cap H)_+)$ we
obtain
$$
\alpha^{N+k+1}(xh_k)\in \alpha(V(O\cap H))\subseteq
O\alpha(O\cap H)\subseteq O\alpha((O\cap H)_+)\;.
$$

Choose $l_{k+1}$ in $(O\cap H)_+$ such that
$\alpha(\alpha^{N+k}(xh_k)l_{k+1}^{-1})\in O$ and put
$h_{k+1}:=h_k\alpha^{-(N+k)}(l_{k+1}^{-1})$.
Then the element $h_{k+1}$ is in $h_0(O\cap H)$, 
and $\alpha^{N+k+1}(xh_{k+1})\in O$.
Further, since $xh_k$ is in the same $H$-coset as $x$,
$\alpha^{N+k+1}(xh_{k+1})\in VH$, by the definition of $N$.

Hence $\alpha^{N+k+1}(xh_{k+1})\in O\cap VH=V(O\cap H)$.
For all natural numbers $j$ less than $k+1$ we have
$\alpha^{N+j}(xh_{k+1})=\alpha^{N+j}(xh_k)\alpha^{j-k}( l_{k+1}^{-1})\in
V(O\cap H)$  as well, showing the existence of our announced sequence.

The sequence $(xh_i)_{i\in\mathbb{N}}$ is bounded, hence has an accumulation
point $xh$ in $xH$.
Then $\alpha^{N+j}(xh)$ is in $V(O\cap H)$ for any natural number $j$ 
because $V(O\cap H)$ is closed. 
In particular, it is in $O$, showing the first claim. 

By the first claim, which we already proved, 
$\{\alpha^n(xh)\colon n\in\mathbb{N}\}$ is bounded. 
Continuity of $\alpha$ and $\alpha^{-1}$ imply that 
the set of accumulation points of $(\alpha^n(xh))_{n\in\mathbb{N}}$
is an $\alpha$-stable subset of $O$. 
Hence each of these accumulation points belongs to $O_0$. 
On the other hand, $(\alpha^n(x))_{n\in \mathbb{N}}$ hence 
$(\alpha^n(xh))_{n\in \mathbb{N}}$ converges to $e$ modulo $H$. 
This can be reformulated using part~3 of Lemma~\ref{conv.mod.H=acc.in.H} to read that each accumulation point of $(\alpha^n(xh))_{n\in\mathbb{N}}$
belongs to $H$. 
We conclude that each accumulation point of
$(\alpha^n(xh))_{n\in\mathbb{N}}$ belongs to $H\cap O_0$. 
Applying part~3 of lemma~\ref{conv.mod.H=acc.in.H} 
once more we see that the sequence $(\alpha^n(xh))_{n\in\mathbb{N}}$ converges to
$e$ modulo the compact $\alpha$-stable subgroup $H\cap O_0$ of $H$, 
and we have established the second claim. 
The lemma is proved. \qed

Returning to the proof of the theorem, assume that the induction 
hypothesis has been established for $i$. 
We apply the lemma with $O^{(i)}_0\cap H$ in place of $H$, 
$O^{(i+1)}$ in place of $O$ and 
$y_i$ in place of $x$.  
We deduce that there is an element $h_i\in O^{(i)}_0\cap H$ 
and an integer $N_{i+1}$ such that 
$\alpha^n(y_ih_i)$ is contained in $O^{(i+1)}$ for each
$n\ge N_{i+1}$.
Furthermore the sequence $\left(\alpha^n(y_ih_i)\right)_{n\in \mathbb{N}}$ converges to $e$ modulo $O^{(i+1)}_0 \cap H$. 
Putting $y_{i+1}:=y_ih_i$, this gives the induction statement for $i+1$ 
proving that the statement holds for all positive integers. 

Since 
$y_{i+1}\in y_i(O_0^{(i)}\cap H)$ and $O^{(i+1)}\leqslant O^{(i)}$, 
$\bigl(y_i(O_0^{(i)}\cap H)\bigr)_{i=1}^\infty$ 
is a decreasing sequence of compact sets 
and $\bigcap_{i=1}^\infty y_i(O_0^{(i)}\cap H)$ is a single point 
because $\bigcap_{i=1}^\infty O^{(i)}$ is trivial. 
Let $\{y\}=\bigcap_{i=1}^\infty y_i(O_0^{(i)}\cap H)$. 
Then for every $n$ 
$$
\alpha^n(y)\in 
\alpha^n\bigl(y_i(O_0^{(i)}\cap H)\bigr)= 
\alpha^n(y_i)\alpha^n(O_0^{(i)}\cap H)
\,.
$$
Since 
$O_0^{(i)}\cap H$ is $\alpha$-stable, 
this set equals 
$\alpha^n(y_i)(O_0^{(i)}\cap H)$, 
which is contained in $O_0^{(i)}$ for every $n\ge N_i$. 
Hence $\alpha^n(y)$ converges to $e$ 
as $n$ tends to  $\infty$. 
Since $y_0=x$, 
$y_{i+1}\in y_iH$ and $H$ is closed, 
$x^{-1}y\in H$ and so the proof is completed by setting $h:=x^{-1}y$.  
\qed

Groups of the form $U_{\alpha/H}$ with $H$ compact arise naturally when
studying (semi-)stable convolution semigroups of probability measures
on the ambient group, see
\cite{Aut(Lie):contr/cp<+Appl(semist*semiG),contrG+semist-meas(p-adicLie)}.
Theorem~\ref{gen_quotient-U} above generalizes Theorem 2.4
in~\cite{contrG+semist-meas(p-adicLie)}. 
Theorem 3.1 in \cite{Aut(Lie):contr/cp<+Appl(semist*semiG)}
covers the case where $H$ is a compact subgroup of a Lie group. 
We pose the question whether Theorem~\ref{gen_quotient-U} 
can be generalised to include the non-metric case. 


The analogue of the above 
Theorem for parabolic groups and their
Levi factors does not hold.
Indeed, there are discrete counterexamples.

\begin{beispiel}
%
Let $G$ be a finitely generated discrete group such that
$G\neq [G,G]=Z(G)$ and $G/[G,G]$ is torsion free.
For example, take $G$ to be the group of integral strict upper
triangular matrices of rank 3.
Take $L:=[G,G]$ and let $g$ be some element of $G$.
Since $G/L$ is abelian, the parabolic group (and its Levi factor) attached
to 
$gL\in G/L$ is the whole group.
We will show that $P_g$ surjects onto $P_{gL}$ only in the trivial
case $g\in L$, providing the desired counterexample.

The assumption that $P_g$ surjects onto $P_{gL}=G$ implies that $G=P_gL=P_gZ(G)=P_g$. 
Hence every element of $G$ has only a finite number, 
$m(x)$ say, of $\langle g\rangle$-conjugates and we 
infer that $[g^{m(x)},x]=e$. 
But $G$ is finitely generated, thus there is a positive $M$ such
that $[g^M,x]=e$ holds for all $x$ in $G$.
In other words, $g^M$ is in $Z(G)$.
Since $G/Z(G)=G/[G,G]$ is torsion free, this implies that $g$ belongs
to the center of $G$, that is, it belongs to $L$ as claimed.
\end{beispiel}

Though Proposition~\ref{quotient-U} does not generalize to parabolic
groups and Levi factors, the following weaker result is obvious.

\begin{satz}\label{quotient-P,M}
Let $p\colon G\to \overline{G}$ be a perfect homomorphism of locally
compact  
groups and let $\alpha$ be an  automorphism of $G$
leaving the kernel of $p$
stable and thus inducing an  automorphism
$\overline{\alpha}$ of $\overline{G}$.
Then $p^{-1}(P_{\overline{\alpha}})= P_\alpha$ and
$p^{-1}(M_{\overline{\alpha}})= M_\alpha$. 
In particular $p(P_\alpha)= P_{\overline{\alpha}}$ and
$p(M_\alpha)= M_{\overline{\alpha}}$
\end{satz}

We conclude the subsection with 
some 
examples illustrating these
concepts.

\begin{beispiel}\label{1stEx:U,P}
The following examples present different types of
behavior which can occur.
The first and second of these should provide the reader
with 
geometrical intuition
on what is going on.\\[-1ex]

\noindent(1)\quad 
Let $k$ be a locally compact, totally disconnected field,
{\it  e.g.\/} the $p$-adic numbers\/ $\mathbb{Q}_p$. 
Let $G$ be
$SL_n(k)$, equipped with the subspace topology in $k^{n^2}$.
Then $G$ is a totally disconnected locally compact group.

In the case where $G=SL_n(\mathbb{Q}_p)$, it follows from
Lemma \ref{red71.2} that, when
computing $U_g$ and $P_g$, we may suppose that $g$ is
semisimple.
That means that, after conjugation, we may assume
that $g$ is diagonal over some finite extension field of\/
$\mathbb{Q}_p$ (also see
Proposition~\ref{scale(redG)}).
We may further assume, that
the valuations of the diagonal entries are in decreasing order.
It is easy to compute $U_g$, $P_g$ and $M_g$ with this
normalization and one finds that $U_g$ is a closed (algebraic)
normal subgroup of $P_g$ consisting
of unipotent matrices, and that $P_g$ is the semidirect product
of $M_g$ and $U_g$. 
If $U_g$ is bounded, we have $U_g=1$.
If all eigenvalues have distinct absolute
value, $M_g$, $P_g$ and $U_g$ will be the groups of diagonal,
upper and strictly upper triangular matrices respectively.
If $\G.$ is any semisimple group, then essentially the same
results
hold for the group $\G k.$: $P_g$ is the group
of rational points of a $k$-parabolic subgroup; $U_g$ is
the group of rational points of its unipotent radical; and
$M_g$ is the group of rational points of the centralizer of
the unique $k$-split torus contained in $P_g\cap P_{g^{-1}}$.
(This is essentially contained in Lemma 2 of \cite{elem:BTR+T}.
One should
note, that the hypothesis that $\G.$ is almost $k$-simple
is not used in its proof 
and that the hypothesis
on the eigenvalues of $\text{Ad}(g)$ is only needed to ensure
that $P_g$ is a proper subgroup.)

The scale function was computed for 
general and special linear groups over local skew fields 
(and some other linear groups) 
by Gl\"ockner in
\cite{scaleF(Ln(local-skew-fields))}. 
This list includes the group $SL_n(\mathbb{Q}_p)$ 
discussed above. 
In that group the subgroup
$\id +M_n(p\mathbb{Z}_p)$, is shown to be tidy for diagonal $g$.
It turns out, that $s(g)$ is the product of the absolute values of
those eigenvalues of $\text{Ad}(g)$ which have absolute value
greater than or equal to $1$.
The scale function is computed for connected semisimple algebraic 
groups over arbitrary
non-Archimedean fields in Proposition~\ref{scale(redG)}.
If the characteristic of the field is $0$ one can alternatively use
Lie methods to compute the scale function as done in
\cite{scaleF(p-adic.Lie)}.\\[-1ex] 

\noindent(2)\quad 
Let $T$ be the homogeneous tree of degree $q+1$.
Taking fixators of finite sets of vertices as basic neighborhoods
of the identity induces a \tdlcG\ topology on $\text{Aut}(T)$.
Each automorphism of $T$ either has a fixed point, not
necessarily
a vertex (elliptic case) or a stable line, which is unique and
is called the axis of $g$ (hyperbolic case).

An elliptic element $g$ is topologically periodic, hence has
$U_g=1$, $P_g=\text{Aut}(T)=M_g$ and trivial scale. 
The stabilizer of any point fixed by $g$ is a tidy subgroup            
since it contains $g$. 
If $g$ is hyperbolic, then it is easily shown that
$P_g$ is the stabilizer of the repelling end $\epsilon_-$
of $g$. 
An automorphism $x$ is in $U_g$ if and only if  for each
$r>0$
there is a point $p(x,r)$
on the axis of $g$ such that all points within distance $r$
of the ray $[p(x,r),\epsilon_-)$ are fixed by $x$.
In this case, $U_g$ is unbounded but not closed (the closure of
$U_g$ is the set of elliptic elements fixing $\epsilon_-$).
The group $M_g$ being the fixator of the two ends fixed by $g$,
we easily see that $M_g\cap U_g$ is nontrivial.
The scale of a hyperbolic element is $q^{l(g)}$, where
$l(g)$ is the
length of the translation induced by $g$ on its axis.
The fixator of a segment of length at least one on the axis
of $g$ is tidy for $g$ (one can show using Lemma
\ref{tidy iff.U0}(3), that
every tidy subgroup for $g$ is essentially of that form).
This example was treated in detail in section 3 of
\cite{tdlcG.structure}.\\[-1ex]

\noindent(3)\quad 
Let $H$ be a \tdlcG\ and $O$ a compact, open subgroup.
For example we may take both $H$ and $O$ equal to a finite group
$F$ carrying the discrete topology.
The shift map $i\mapsto i+1$ on $\mathbb{Z}$
induces an automorphism $\sigma$ of the restricted
product $G:=\prod_{i\in\mathbb{Z}} H|O$.
Any compact subset stable under the shift $\sigma$ is contained
in
$O^\mathbb{Z}=\prod_{i\in\mathbb{Z}} O$ and any open stable
set contains $O^\mathbb{Z}$.
Therefore, $O^\mathbb{Z}$ is the unique
$\sigma$-stable compact, open subgroup of $G$, 
and thus the only subgroup tidy for $\sigma$.
As a corollary, $\sigma$ has scale $1$. 

The parabolic subgroup $P_\sigma$ and Levi-factor $M_\sigma$
are both equal to $O^\mathbb{Z}$. 
The contraction group $U_\sigma$ 
is the subgroup of all
$(x_i)_{i\in\mathbb{Z}}\in O^\mathbb{Z}$ 
such that
$x_i\too_{i\to-\infty}e$.
This is a nontrivial and bounded group which is not closed. 
\end{beispiel}

\subsection{The difference between shrinking and contracting}
\label{subsec:difference}

From now on we assume that the ambient group $G$ is totally disconnected
and \emph{metrizable\/} unless explicitly stated otherwise.
We establish the connection between that part, $V_-$,
of a compact, open subgroup $V$ on which $\alpha$ is shrinking and the
corresponding contraction group $U_\alpha$.  
The notion of contractivity is stronger than that of  shrinking, as the first lemma demonstrates.

\begin{lemma}\label{U<V-}
Let $V$ be a compact, open subgroup of a \tdlcG. 
Then $U_\alpha\leqslant V_{--}$.
\end{lemma}
\proof{} Take $v\in U_\alpha$. 
Then there is an $N$ such that
$\alpha^{n}(v)$ belongs to (the $e$-neighborhood) $V$ for every $n\geq N$.
Thus $\alpha^{N}(v)\in V_-$ and $v\in V_{--}$. \qed

Since $V_-\subseteq P_\alpha$ (see Remark \ref{V- in P}),  the following
result is an immediate consequence of the lemma and Proposition
\ref{U normal P}.

\begin{corollar}
Let $V$ be a compact, open subgroup of a \tdlcG. 
Then 
$U_\alpha\normaleq V_{--}$,
$U_\alpha \cap V_{-}\normaleq V_{-}$ and $U_\alpha \cap V_{0}\normaleq V_{0}$.
\qed
\end{corollar}

We now make the interdependence between $U_\alpha$ and $V_{--}$
more precise. 
Clearly, $U_\alpha V_0\leqslant V_{--}$ and
Lemma~\ref{cp-contr} implies that $V_{--}\subseteq U_\alpha$
only when $V_0=1$.
In fact, $V_0$ is the `difference' between $V_{--}$ and $U_\alpha$.

\begin{satz}\label{V--=}
For a compact, open subgroup $V$, 
$V_{--}=U_\alpha V_0$.
\end{satz}
\proof{} 
We are left to show $V_{--}\subseteq U_\alpha V_0$. 
For each $n\in\mathbb{N}$ and $v\in V$ 
we have 
$\alpha^n(v)\in \bigcap_{j= -N}^N \alpha^j(V)$ 
for all sufficiently large $n$. 
Since $\bigcap_{j= -N}^N \alpha^j(V)$ is compact 
and decreases to $V_0$ with $N$, 
$(\alpha^n(v))_{n\in \mathbb{N}}$ 
converges to $e$ modulo $V_0$ 
for all $v\in V_{--}.$ 
The result follows from Theorem~\ref{gen_quotient-U}.  
\qed 

The next result may be considered a global version of the above 
propositon. 

\begin{corollar}\label{P=MU}
$M_\alpha U_\alpha= P_\alpha$
, in other words $U_{\alpha/M_\alpha}=P_\alpha$. 
\end{corollar}
\proof{}
Using Proposition~\ref{U normal P} it suffices to show 
$P_\alpha\subseteq U_{\alpha/M_\alpha}$. 
Let $v$ be in $P_\alpha$. 
By definition of $P_\alpha$ the sequence
$(\alpha^n(v))_{ n\in \mathbb{N}}$ is bounded. 
Let $w$ be an accumulation
point, i.e. a limit of a subsequence.
Then for any $m\in \mathbb{Z}$ the point $\alpha^m(w)$ as well is
a limit of a subsequence and belongs to the compact set
$\ol{\{\ \alpha^n(v)\colon n\in \mathbb{N}\}}$. 
This means $w\in M_\alpha$. 
Since the set $\{\alpha^n(v)\colon n\in\mathbb{N}\}$ is bounded and 
$M_\alpha$ is closed,  
we may apply part~3 of Lemma~\ref{conv.mod.H=acc.in.H} to conclude 
that $v\in U_{\alpha/M_\alpha}$ as claimed. 
\qed

As might be expected, the topology of the quotient space
$P_\alpha/\overline{U}_\alpha$ is induced from $M_\alpha$.

\begin{lemma}\label{M/U0}
The natural homomorphism $M_\alpha\to P_\alpha/\overline{U}_\alpha$
is an identification.
\end{lemma}
\proof{}
The homomorphism is continuous and surjective.
To show that is open, it suffices to find an open subgroup of $M_\alpha$
such that the restriction of the homomorphism to this subgroup is open.

Let $V$ be a subgroup tidy for $\alpha$ in the ambient group.
Then $V_0$ equals $V\cap M_\alpha$ and therefore is a compact open subgroup
of $M_\alpha$. 
We first claim that the image $O$ of $V_0$ is open:
The inverse image of $O$ is
$V_0\overline{U}_\alpha\leqslant P_\alpha$.
The equations $V_0 U_\alpha=V_{--}\geqslant V_- =V\cap P_\alpha$ show
that $V_0U_\alpha$ is an open hence closed subgroup of $P_\alpha$.
It follows that $V_0\overline{U}_\alpha$ must equal $V_0U_\alpha$ and
is therefore open. 
This proves that $O$ is open as claimed.

By the open mapping theorem 
(\cite[II.(5.29)]{GMW115})
the map $V_0\to O$ is necessarily open. 
Since $O$ was shown to be open, this implies that the map 
$V_0\to P_\alpha/\overline{U}_\alpha$
is open. 
\qed

Lemma \ref{U0:<>UPM} and Corollary \ref{U0:=} will show that
the kernel of the homomorphism 
$M_\alpha\to P_\alpha/\overline{U}_\alpha$ is compact. 
It then follows that this homomorphism is a perfect map.

\subsection{Reinterpretation of the scale function}
We will investigate the links between contraction groups,
parabolics and tidy subgroups further after giving alternative
descriptions of the scale of an automorphism. 

The next two lemmas are immediate consequences of Example \ref{V=V-} and
Remark \ref{V- in P}.

\begin{lemma}\label{V-H}
Let $V$ be a compact, open subgroup of a \tdlcG\   
tidy for the automorphism $\alpha$.
Then for all closed $\alpha$-stable subgroups $H$ of $P_\alpha$
$$
(V\cap H)_- = V_-\cap H =V\cap P_\alpha\cap H =V\cap H
$$
is tidy for $\alpha$ in $H$. \qed
\end{lemma}

\begin{lemma}\label{V-N}
Let $V=V_-$ be a compact, open subgroup of the \tdlcG\ $H$ and let
$N\normaleq H$ be stable under the automorphism $\alpha$ of $H$.
Then the image $q(V)\subseteq H/N$ under the canonical map satisfies
$$
q(V_-)\subseteq q(V)_-\subseteq q(V)=q(V_-)\ ,
$$
and hence is tidy for the induced automorphism
$\ol{\alpha}\colon H/N\to H/N$. \qed
\end{lemma}

The reason for the name 'scale' 
is that $s(\alpha)$ is the factor by which $\alpha$ 
scales up $V_+$ when $V$ is tidy. 
Thus $s(\alpha)$ is just the value of the modular function 
at the restriction of $\alpha$ to $V_{++}$. 
The next result extends this interpretation of the scale.

\begin{satz}\label{s=mod()}
Let $N{\normaleq} H$ be 
a $\alpha$-stable closed subgroups of $P_\alpha$ and
let $V$ be tidy for $\alpha$ in the ambient group $G$.
Then writing $q$ for the canonical map $H\to H/N$ and $\ol{\alpha}$,
respectively $\alpha_|$, for the induced automorphisms on $H/N$ and $N$,
we have
\begin{itemize}
\item[(1)] $s_H(\alpha^{-1})=\Delta_H(\alpha^{-1})$
\item[(2)] $s_H(\alpha^{-1})=s_{H/N}(\ol{\alpha}^{-1})s_N(\alpha^{-1}_|)$
\item[(3)] $s_G(\alpha^{-1})=s_{P_\alpha}(\alpha^{-1})=s_{V_{--}}(\alpha^{-1})=
       s_{\ol{U}_\alpha}(\alpha^{-1})\;.$
\end{itemize}
\end{satz}
\proof{} 
``(1)'': 
By Lemma \ref{V-H} we have $(V\cap H)_- = V\cap H$
and this group is tidy for $\alpha$ in $H$. 
Hence $\Delta_H(\alpha^{-1})=s_H(\alpha^{-1})$. 

\noindent{``(2)'':} 
By Lemma \ref{V-N} we have $q(V\cap H)_- = q(V\cap H)$
and this group is tidy for $\ol{\alpha}$ in $H/N$. 
Hence
$\Delta_{H/N}(\ol{\alpha}^{\;-1})=s_{H/N}(\ol{\alpha}^{\;-1})$.
Additionally,  
applying (1) with $H:=N$ 
implies
$\Delta_{N}(\alpha_|^{-1})=s_{N}(\alpha_|^{-1})$.

Since the modular function satisfies an equation of the form to be proved, 
substituting the values of the scale function for those of the modular function 
in that equation proves (2). 

%

\noindent{``(3)'':} 
We have $s_G(\alpha^{-1})=|\alpha^{-1}(V_-):V_-|=
s_{P_\alpha}(\alpha^{-1})$ since $V_-=V\cap P_\alpha=(V\cap P_\alpha)_-$.
The same argument works when $P_\alpha$ is replaced by $V_{--}$.
Now (2) applied  with $H:=P_\alpha$ and $N:=\ol{U}_\alpha$ gives
$s_{P_\alpha}(\alpha^{-1})=s_{P_\alpha/\ol{U}_\alpha}(\ol{\alpha}^{\,-1})
s_{\ol{U}_\alpha}(\alpha^{-1}_|)$, leaving us to prove
$s_{P_\alpha/\ol{U}_\alpha}(\ol{\alpha}^{\,-1})=1$. 

As 
the product of the scale of the restriction 
of an automorphism to a normal subgroup 
and 
the scale of the induced automorphism on the quotient 
always 
divides 
the scale of the automorphism  
by Proposition 4.7 in \cite{furtherP(s(tdG))},  
it is enough to show
$s_{M_\alpha}(\alpha^{-1})=1$, since $P_\alpha/\ol{U}_\alpha$ is a
quotient of $M_\alpha$ by Lemma \ref{M/U0}. 
We compute
$$
s_{M_\alpha}(\alpha^{-1})=|\alpha^{-1}(M_\alpha\cap V)_-:(M_\alpha\cap V)_-|
$$
Observing that
$$
M_\alpha\cap V=
P_\alpha\cap P_{\alpha^{-1}}\cap V=V_+\cap V_-=V_0,
$$
this implies $s_{M_\alpha}(\alpha^{-1})=|V_0:V_0|=1$, as had to be shown.
\qed

~\\
Combination of (3) and (1) above implies that
$s_G(g)=\Delta_{\overline{U}_{g^{-1}}}(g)$.
This enables us to compute the scale function of the group of
rational points of a semisimple algebraic group $\G.$
over a local field of positive characteristic. 
We start with the following lemma.

\begin{lemma}
\label{lem:commuting}
If $\alpha$ and $\beta$ are two commuting automorphisms
then
\begin{itemize}
\item[(1)]
$s(\alpha\beta)\le s(\alpha)s(\beta)$,
\item[(2)]  $s(\beta)=1=s(\beta^{-1})$ implies
$s(\alpha\beta)=s(\alpha)$.
\end{itemize}
\end{lemma}
\proof{}
As (2) follows from (1), it suffices to prove (1). 
We may suppose that $\alpha$ and $\beta$ are 
inner. 
Using Theorem 3.4 in \cite{A(JordanF(tdlcG))+simulT}, 
choose a compact open subgroup $V$, 
which is tidy for $\alpha$ and $\beta$. 
The result follows from Proposition~\hbox{A.2}
in \cite{top(HeckeP+HeckeC*)}.
\qed

We now treat a slightly more general case than announced in
Example~\ref{1stEx:U,P}.

\begin{satz}\label{scale(redG)}
Let $k$ be a nonarchimedean local field and let $\mathbf{G}$ be a
Zariski-connected reductive $k$-group. 
For any element $g$
of $\mathbf{G}(k)$ its scale $s(g)$
equals the product of the absolute values of those eigenvalues
of $Ad(g)$, whose
valuation is greater than $1$ (counted with their proper multiplicities).
\end{satz}
\proof{}
Observe first that the claim will be true for an element $g$ whenever
it is true for a positive power of $g$. 
Likewise, using part (2) of
Lemma \ref{lem:commuting}, we see that we may replace $g$ by $g'$
whenever $g^{-1}g'$ is a compact element commuting with $g$. 
If $k$ has positive characteristic some positive power of $g$ is semisimple. 
If the characteristic of $k$ is $0$,
$g$ has a Jordan-Chevalley decomposition $g=g_sg_u=g_ug_s$ with $g_s$, $g_u\in \mathbf{G}(k)$
where $g_s$ is semisimple and $g_u$ is unipotent. 
Since unipotents are
compact elements, we may replace $g$ by its semisimple part $g_s$. 
So, regardless of the
characteristic, we may assume from the outset  that $g$ is semisimple.

Let $\mathbf{S}$ be the smallest torus of $\mathbf{G}$ containing $g$, that
is, let $\mathbf{S}$ be the Zariski-closure of the group generated by $g$.
From this characterization it is immediate that $\mathbf{S}(k)$ is
Zariski-dense in $\mathbf{S}$. 
By the Galois-criterion $\mathbf{S}$ is defined over $k$. 
Its largest split, respectively anisotropic, tori $\mathbf{S}_d$ and
$\mathbf{S}_a$ are defined over $k$ as well.

Write $g$ as $g'a$, where $g'$ is in $\mathbf{S}_d(k)$ and $a$ is in
$\mathbf{S}_a(k)$. 
By our introductory remarks we may suppose that $g=g'$ is
in fact
contained in a $k$-split torus.

We are going to use the Remarque after Corollaire 3.18 in \cite{ihes27}.
Let $\sigma$ be the  set of roots of $\Phi(\mathbf{S}_d,\mathbf{G})$ with
$|\sigma(g)|>1$. 
It is a closed (even connected) subset and by
\cite[Chapitre VI,\jura 1, Proposition 22]{gaLie4-6} there is an ordering on $\Phi(\mathbf{S}_d,\mathbf{G})$ such that all the roots in $\sigma$ are
positive. 
Hence $\sigma$ is unipotent. 
The set of $k$-rational points of the associated unipotent
group $\mathbf{U}_\sigma$ is equal to $U_{g^{-1}}$, hence $U_{g^{-1}}$
is already closed.
This follows from the fact that $\text{Lie}(\mathbf{G})$ is the direct sum of
$\text{Lie}(\mathbf{U}_\sigma)$, 
$\text{Lie}(\mathbf{U}_{-\sigma})$ and
the sum of the eigenspaces of $g$ to eigenvalues with valuation $1$.  
In particular each  eigenspace of $g$ to an
eigenvalue with valuation greater than $1$ is contained in
$\text{Lie}(\mathbf{U}_\sigma)$.

Arrange the elements $b_1,\ldots,b_m$ of $\sigma$ in increasing order and put $\sigma_i:=\{b_i,\ldots,b_m\}$. 
The group $\mathbf{U}_\sigma$ is defined and
split over $k$. 
More precisely,
$\mathbf{U}_\sigma=\mathbf{U}_{\sigma_1}\geqslant \cdots \geqslant
\mathbf{U}_{\sigma_m}$ is a filtration by unipotent
$k$-groups whose successive quotients admit a structure of $k$-vector
space such that any element $x$ of $\mathbf{S}_d$ acts by scalar multiplication with
$b_i(x)$ on $\mathbf{U}_{\sigma_i}/\mathbf{U}_{\sigma_{i+1}}$. 
It follows that $\Delta_{U_{g^{-1}}}(g)= \Delta_{\mathbf{U}_\sigma(k)}(g)=\prod_{i=1}^m
b_i(g)^{d_i}$, where $d_i=\dim(\mathbf{U}_{\sigma_i}(k)/
\mathbf{U}_{\sigma_{i+1}}(k))$. 
This is exactly what we claimed. \qed

As a further result, contraction groups are usually unbounded.

\begin{satz}\label{U bounded iff}
The following statements are equivalent:\\
\begin{tabular}{rlrl}
{(1)} & $s_G(\alpha^{-1})=1$.&&\\[.7ex]
{(2)} & $U_\alpha$ is bounded.&&\\[.7ex]
{(3)} & $\forall O$ tidy for $\alpha$ : $O\subseteq P_{\alpha^{-1}}$.
    & \ {(A)} & $U_\alpha\subseteq M_\alpha$.\\[.7ex]
{(4)} & $\exists O$ tidy for $\alpha$ : $O\subseteq P_{\alpha^{-1}}$.
    & \ {(B)} & $M_\alpha=P_\alpha$.\\[.7ex]
{(5)} & $P_{\alpha^{-1}}$ is open.
    & \ {(C)} & $\ol{U}_\alpha=U_0:=
      \ol{U}_\alpha\cap\ol{U}_{\alpha^{-1}}$.\\
\end{tabular}
\end{satz}
\proof{}
We will first prove equivalence of the statements in the
left column. 
The proof that the characterizations in the
right column are equivalent uses some results yet to come but are listed
here for convenience of reference. 
These equivalences will not be used until later.

\noindent``(2) $\Rightarrow$ (1)'': 
This is immediate combining
parts (1) and (3) of Proposition \ref{s=mod()} (with
$H:=\ol{U}_\alpha$).
We prove that (1) implies (2) and (3).\\
Take $O$ tidy for $\alpha$.
The equality $1=s_G(\alpha^{-1})=|\alpha^{-1}(O):\alpha^{-1}(O)\cap O|$
is equivalent
to $O\supseteq \alpha^{-1}(O)$ and to $O_+=O$.
The first of these implies
$O_{--}\subseteq \bigcup_{i\in \mathbb{N}}\alpha^{-i}(O)\subseteq O$ hence
that $U_\alpha\subseteq O$ is bounded, giving (2).
Further $O=O_+=O\cap P_{\alpha^{-1}}$ gives $O\subseteq P_{\alpha^{-1}}$,
hence (3) follows.

Obviously $(3)\Rightarrow (4)\Rightarrow (5)$.
Finally to prove that (5) implies (1) take $O$ tidy for $\alpha$ in
$P_{\alpha^{-1}}$.
Since we assume that $P_{\alpha^{-1}}$ is open, $O$ is tidy for
$\alpha$ in the ambient group $G$ as well and we get
$O_+=O\cap P_{\alpha^{-1}}=O$ which is
equivalent to $1=s_G(\alpha^{-1})$ as has already been seen.

We now prove equivalence of (2) and the statements in the right column.
If we assume (2) then using $\alpha$-invariance of $U_\alpha$
we get that $\alpha^\mathbb{Z}(U_\alpha)=U_\alpha$ is bounded.
The definition of $M_\alpha$ then gives $U_\alpha\subseteq M_\alpha$,
that is (A).

The statements (A) and (B) are equivalent thanks to Corollary \ref{P=MU}.
For the remaining implications we need some results yet to be proved.
If we assume (B) then since $M_\alpha$ is closed
$\ol{U}_\alpha=\ol{U}_\alpha\cap M_\alpha$ which equals $U_0$ by
Lemma \ref{U0:<>UPM} and we have derived (C).
Assuming (C) we have $U_\alpha$ contained in $U_0$, which is a compact
group by Corollary \ref{U0:=} giving (2).
The proof is complete.
\qed

\subsection{Small tidy subgroups}
The contraction group, $U_\alpha$, is closed if and only if there are
arbitrarily small subgroups tidy for $\alpha$. 
This and numerous other equivalences are established below in Theorem \ref{U0=1 iff}. 
Careful examination of the tidying procedure, Algorithm \ref{algo}, is required for the proofs.

\begin{corollar}[to Proposition~\ref{V--=}] \label{-tidyup=}
If $V$ is a tidy subgroup produced by the Algorithm
\ref{algo} from the compact, open subgroup $O$, then
$V_{--}=\ol{O_{--}}$.
\end{corollar}
\proof{} We go through the steps of the Algorithm \ref{algo}.
Step~1 produces a subgroup $O'$ such that
$O'_-=\alpha^n(O_-)$ hence
$$
O'_{--}=O_{--}\;.
$$

The group $L=\ol{\mathcal{L}}$ produced in step~2 of the algorithm
equals $V_0=O''_0$ by  Lemma 3.8 in \cite{furtherP(s(tdG))} and is
easily seen, by its definition, to be contained in $\ol{O'_{--}}$.
Hence, by Proposition~\ref{V--=},
$$
V_{--} = U_\alpha V_0 = U_\alpha L
\leqslant U_\alpha \ol{O'_{--}} = \ol{O'_{--}}=\ol{O_{--}}\;.
$$

To show the reverse inclusion we proceed to step~3 of the algorithm. 
Lemma 3.4 in \cite{furtherP(s(tdG))} implies
$$
O^*_{--}= O'_{--}\,.
$$
From $O^*\leqslant V$ we get immediately  $O^*_-\leqslant V_-$ and
$O^*_{--}\leqslant V_{--}$ hence
$$
\ol{O_{--}}= \ol{O'_{--}}= \ol{O^*_{--}}\leqslant 
\ol{V_{--}}=V_{--}
$$
since $V$ is tidy. 
We are done. \qed

As an immediate consequence, we have the following characterisation of the
closure of $U_\alpha$.

\begin{theorem}\label{Ualpha bar}
$\ol{U}_\alpha= \bigcap \{V_{--}\colon V\ \text{is tidy for $\alpha$}\, \}\,.$
\end{theorem}
\proof{}
We already know $U_\alpha\subseteq V_{--}$ for every compact open
subgroup $V$. 
Since $V_{--}$ is closed for $V$ tidy, the inclusion
$\ol{U}_\alpha\subseteq \bigcap\{V_{--}\colon V\ 
\text{is tidy for $\alpha$}\, \}$ follows.

Now let $v\notin \ol{U}_\alpha$. 
We have to find a subgroup $V$ tidy for
$\alpha$ such that $v\notin V_{--}$.
There is a compact open subgroup $O$ such that
$vO\cap \ol{U}_\alpha = \emptyset$.
Then $vO\cap \ol{U}_\alpha O_0= \emptyset$ and it follows that
$$
v\notin \ol{U_\alpha O_0} = \ol{O_{--}}.
$$
Thanks to Corollary \ref{-tidyup=}, $v\notin V_{--}$,
where $V$ is the tidy subgroup constructed from $O$.
We are done. \qed

One would expect the set of elements where $\alpha$ and $\alpha^{-1}$
are contracting to be small.

\begin{corollar}\label{U0:=}
The group $U_0$, defined as $\ol{U}_\alpha\cap \ol{U}_{\alpha^{-1}}$, is
equal to
$$
\bigcap \{V_0\colon V\ \text{is tidy for $\alpha$} \},
$$
 and hence to
$$
\bigcap \{V \colon V\ \text{is tidy for $\alpha$} \}
$$
as well.  
In particular, it is compact.
\end{corollar}
\proof{} First observe that if $V$ is tidy, then so is $\alpha^n(V)$ for any
integer $n$.
Therefore
$$
\bigcap \{V \colon V\ \text{is tidy for $\alpha$}\}=
\bigcap_{V \text{tidy}}\bigcap_{n\in \mathbb{Z}} \alpha^n(V)=
\bigcap \{V_0\colon V\ \text{is tidy for $\alpha$}\}.
$$
This is a compact group. 
It remains to show that it equals $U_0$.

For this, note that, by Theorem \ref{Ualpha bar}, $U_0$ is the
intersection of all pairs
$V_{--}$ and
$W_{++}$ where $V$ and $W$ run through all tidy subgroups.

Since $V_0\subseteq V_{--}$ and $W_0\subseteq  W_{++}$, it follows
that
$$
\bigcap\{V_0 \colon V \ \text{is tidy for } \alpha\}\subseteq U_0.
$$
The next lemma implies the reverse inclusion.

\begin{lemma}\label{V++cap V--}
If $V$ be a tidy subgroup of a \tdlcG. 
Then $V_{++}\cap V_{--}=V_0$.
\end{lemma}
\proof{} As already seen $V_0=V_+\cap V_-\subseteq V_{++}\cap V_{--}$ always.
Conversely if $x\in V_{++}\cap V_{--}$ then
$$ 
x\in \alpha^k(V_+)\cap \alpha^l(V_-)\ \text{for all}\
k\ge k_0\ge 0\ge l_0\ge l\ .
$$
Now $\alpha^{-k}(x)\in V_+\subseteq V\ \forall k\ge k_0$ and
 $\alpha^{-l}(x)\in V_-\subseteq V\ \forall l\le l_0$.
Therefore $x\in \mathcal{L}$. 
Assuming $V$ to be tidy gives 
$\mathcal{L}=V_0$, implying that $x\in V_0$.
The proof of the corollary is complete.
\hspace*{\fill}$\square\ \square$

We will now turn to the question of existence of arbitrarily small
tidy subgroups for an automorphism $\alpha$.
Corollary \ref{U0:=} shows that $U_0$ is an obstruction
to their existence. 
We attempt now to make this more precise.

\begin{lemma}\label{U0:<>UPM}
$\ol{U}_\alpha\cap P_{\alpha^{-1}} = \ol{U}_\alpha\cap M_\alpha=
U_0 $.
\end{lemma}
\proof{} The proof is an easy direct computation using
Theorem \ref{Ualpha bar} and that 
$V_0=V\cap M_\alpha$ for tidy $V$:
\begin{eqnarray*}
U_0 &=& \ol{U}_\alpha\cap \ol{U}_{\alpha^{-1}}\subseteq
 \ol{U}_\alpha\cap P_{\alpha^{-1}}= \ol{U}_\alpha\cap M_\alpha
 =\bigcap_{V\,\text{tidy}} V_{--}\cap M_\alpha =\\
&=& \bigcap_{V\,\text{tidy}} (V_{--}\cap M_\alpha)= \bigcap_{V\,\text{tidy}}
\bigl(\bigcup_{i\in \mathbb{N}} \alpha^{-i}(V\cap P_\alpha)\cap M_\alpha\bigr)
=\\
&=&\bigcap_{V\,\text{tidy}}\bigcup_{i\in \mathbb{N}}
  \alpha^{-i}(V\cap M_\alpha) =
\bigcap_{V\,\text{tidy}}\alpha^{-i}(V_0)=\bigcap_{V\,\text{tidy}} V_0=
U_0\,.
\end{eqnarray*}
\qed

As a Corollary we get a nice factorization of $\ol{U}_\alpha$.

\begin{corollar}
\label{cor:Ualphabarfactor}
$\ol{U}_\alpha=U_0 U_\alpha$,  
In other words $\ol{U}_\alpha= U_{\alpha/U_0}$. 
\end{corollar}
\proof{} 
Let $\alpha_1$ be the restriction of $\alpha$ to $\ol{U}_\alpha$. 
Then $U_{\alpha_1}=U_\alpha$, 
$M_{\alpha_1}=U_0$ by Lemma~\ref{U0:<>UPM} 
and $P_{\alpha_1}=\ol{U}_\alpha$ 
and the claim follows by Corollary~\ref{P=MU}. 
\qed

After a further lemma we will be able to characterize the automorphisms
with arbitrarily small tidy subgroups.
Section 2 of \cite{A(JordanF(tdlcG))+simulT} describes a tidying
procedure which differs from Algorithm \ref{algo} by taking in step~3
the group $K:=\ol{U_\alpha\cap P_{\alpha^{-1}}}\subseteq U_0$ instead of
$L$ (constructed in step~2, which becomes superfluous). 
That is, we now put 
$O^*:=\{x\in O'\colon kxk^{-1}\in O'K\,\forall k\in K\}$ 
and $O'':=O^* K$. 
According to Lemma 3.3(1) in \cite{A(JordanF(tdlcG))+simulT},
$O^*$ is a compact open subgroup. 
It is smaller than $O$. 
We work with this algorithm in the proof of the lemma below.

It is proved in {\it  loc.cit.\/}, that the outputs of the two algorithms
are the same.
One may wonder, whether the reason for this is that $K=U_0$.
It follows from (2) and (1) below, that this is indeed the case.

\begin{itemlemma}\label{K=U0}\label{tidy iff.U0}
\begin{itemize}
\item[(1)] $ U_\alpha\cap P_{\alpha^{-1}}= U_\alpha\cap M_\alpha =
           U_\alpha\cap U_0$.
\item[(2)] $\ol{U_\alpha\cap M_\alpha}=U_0$; in particular, $U_0$ is
           the unique tidy subgroup for $\alpha$ restricted to $U_0$.
\item[(3)] A compact, open subgroup is tidy if and only if  it satisfies (T1)
           and contains $U_0$.
\end{itemize}
\end{itemlemma}
\proof{} The first statement is easily derived from the descriptions of
$U_0$ in Lemma
\ref{U0:<>UPM}.

From (1) we get that $\ol{U_\alpha\cap M_\alpha}$ equals
the group
$K$ introduced above. 
It follows from Lemma \ref{U0:<>UPM} that
$K\subseteq U_0$, so assume  we are given an element $x$ not in $K$.
Since $K$ is closed, there  is a compact open subgroup $O$ such that
$xO\cap K$ is empty.  
In other words $x$ does not belong to $KO$.

Applying the modified algorithm to $O$, 
we have $KO\supseteq KO^*=O''$ is a tidy subgroup, and $x\notin KO$ implies
that $x$ is not in $U_0$ by Corollary \ref{U0:=}.
Thus $U_0\subseteq K$ and therefore $U_0=K=\ol{U_\alpha\cap M_\alpha}$
as claimed in (2).

To see the remaining part of sub-claim (2) we compute the group $U_0$
inside $U_0$:
it is $\ol{U_\alpha\cap U_0}\cap \ol{U_{\alpha^{-1}}\cap U_0}$
by its definition, and this expression
equals $\ol{U_\alpha\cap M_\alpha}\cap \ol{U_{\alpha^{-1}}\cap M_\alpha}=
U_0\cap U_0=U_0$ by what is already known.
Corollary \ref{U0:=} implies then, that there can be no smaller
tidy subgroups for the restriction of $\alpha$ to $U_0$ than $U_0$
itself.

We turn to (3).
The conditions are clearly necessary for a subgroup $O$ to be tidy.
If they are satisfied, the variant of the tidying-up
procedure described above leaves the compact open subgroup unchanged.
This shows that $O$ is tidy. \qed

In all examples examined so far, $U_0$ is equal to
$\overline{U_\alpha\cap U_{\alpha^{-1}}}$.
We do not know whether this holds in general.

The existence of arbitrarily small tidy subgroups for
a given automorphism has many equivalent reformulations.

\begin{theorem}\label{U0=1 iff}
The following conditions are equivalent:\\
$(1)$ $U_\alpha$ is closed.\hspace{5pt}\hspace{2.67em}
$(1')$ $U_{\alpha^{-1}}$ is closed.\\[.7ex]
$(2)$ $U_\alpha\cap M_\alpha=1$.\hspace{2.67em}
$(2')$ $U_{\alpha^{-1}}\cap M_\alpha=1$.\\[.7ex]
$(3)$ $\ol{U}_\alpha\cap M_\alpha=1$.\hspace{2.67em}
$(3')$ $\ol{U}_{\alpha^{-1}}\cap M_\alpha=1$.\\[.7ex]
$(4)$ $U_0=1$.\hfill
$(4')$ There are arbitrarily small tidy subgroups for $\alpha$.\\[.7ex]
$(5)$ For all compact, open $O\leqslant G$ there is a $k$ such that
$^kO$ is tidy.\\[.7ex]
$(5')$ A compact, open subgroup satisfying (T1) satisfies (T2) as well.\\[.7ex]
$(6)$ $P_\alpha=M_\alpha\ltimes U_\alpha$ topologically.\hfill
$(6')$ $P_{\alpha^{-1}}=M_\alpha\ltimes U_{\alpha^{-1}}$ topologically.\\
\end{theorem}
\proof{} ``$(5')\iff (5)$'':\
Assume $(5')$. 
Run step~1 of the algorithm \ref{algo} on the given $O$
to obtain another compact, open subgroup of the stated form satisfying
(T1) hence (T2) by assumption.
Conversely assume that $O$ satisfies (T1). 
By assumption $O':= \bigcap_{i=0}^k \alpha^i(O)$ is tidy.
As already used several times  $O'_{+}=O_+$ and $O'_-=\alpha^k(O_-)$ hence
$O_{++}=O'_{++}$ and $O_{--}=O'_{--}$ are closed since $O'$ is tidy.
Hence $O$ is tidy and (T1) is shown to imply (T2).

Since the procedure applied to $O$ in $(5)$ shrinks the compact, open
subgroups in question, we have that $(5)\Rightarrow (4')$, which evidently
implies $(4)$.
Lemma \ref{U0:<>UPM} shows, that $(4)$, $(3)$ and $(3')$ are all
equivalent.

If we assume that $U_\alpha$ is not closed,
then $1\neq U_0$ by Corollary \ref{cor:Ualphabarfactor}. 
Hence $M_\alpha \cap \ol{U}_\alpha\ne 1$ and we have shown that $(3)$ implies $(1)$. 
By symmetry, $(3')$ implies $(1')$ as well.

Assuming  $(1)$ we find for any compact, open subgroup $V$ that
$V_{--}=V_0 U_\alpha$ is closed, since $V_0$ is compact.
Assuming $V$ satisfies (T1) we get from \cite[Lemma 3(b)]{tdlcG.structure}
that $V_{++}$ is closed as well. 
This means that property (T2) is automatic once (T1) is known to hold, 
i.e. $(5')$. 
By symmetry again $(1')$ implies $(5')$.
So far we have established equivalence of all conditions
listed except $(2)$, $(2')$, $(6)$ and $(6')$.

It is trivial that $(3)$ implies $(2)$ and $(3')$ implies $(2')$.
Lemma \ref{K=U0}$(2)$ gives that $(2)$ implies $(4)$.
Symmetrically $(2')$ implies $(4)$. 
Together with $(4)\iff (3)$,
this leaves to prove the equivalence of $(6)$ and $(6')$ with the other
conditions:

By symmetry it suffices to prove that (6) is equivalent to (2).
Property $(6)$ evidently implies $(2)$.
Assume $(2)$ to conclude that $P_\alpha$ is the semidirect product of
$M_\alpha$ and $U_\alpha$ as an abstract group. 
We may then apply Proposition~6.17 from \cite{unif-str(tG+*/)}. 
Our claim is that (6) holds, which is statement (d) in that Proposition. 
It is shown there that it is equivalent to the statement (b) that the map $M_\alpha\to P_\alpha/U_\alpha$
obtained by restriction of the quotient map modulo the normal subgroup
$U_\alpha$ is a topological isomorphism.
Since we know that (2) implies~(1), this map is the one considered in
Lemma~\ref{M/U0}, where we showed that it is an identification in general.
Under the assumption of (2) however, the kernel of this map is trivial and
(b), thus (6), is established. 
We are done.
\qed

We now list some examples 
where the conditions of Theorem~\ref{U0=1 iff} hold.

\begin{bemerkung}\label{Ex(closedU)}
All contraction groups for general/inner automorphisms are closed in the
following cases.

\noindent(1)\quad  
Groups with trivial contraction groups, among them
\begin{enumerate}
\item[(a)] 
discrete groups, with respect to \textbf{all} and
\item[(b)]
SIN-groups with respect to \textbf{inner} automorphisms.
\item[(c)]  
MAP-groups with respect to \textbf{inner}
automorphisms.
(%
The von-Neumann-kernel of a locally compact group
contains all contraction groups of
inner automorphisms:
Let $\rho$ be a continuous finite-dimensional unitary
representation of the group. 
Consider an element, $x$, say.
Since $\rho(x)$ is unitary, the corresponding eigenvalues
have absolute value $1$.
A trivial modification of Lemma II.(3.2)
from \cite{EdM3.17} shows, that $U_x$ fixes each eigenvector
of $\rho(x)$. 
Since there is a basis of the representation space
consisting of eigenvectors of $\rho(x)$, we are done.)
\item[(d)]
nilpotent-by-compact groups. 
These are distal by
\cite[Proposition 3]{distalP(G)+growth}. 
This, by its definition
(loc cit, introduction) easily implies that all contraction groups
of \textbf{inner} automorphisms are trivial.
\end{enumerate}
\noindent(2)\quad
Any (totally disconnected) locally compact group having an
      open subgroup satisfying the ascending chain condition on its
      closed subgroups. 
      Thanks to (\cite[Lemma 3.2]{Mautner(p-adic.Lie)})
      the group then satisfies the condition (1) of Theorem \ref{U0=1 iff}
      with
      respect to \textbf{any} automorphism
      (\cite[Lemma 3.2]{Mautner(p-adic.Lie)}). 
       This criterion applies to
      any $p$-adic Lie group (\cite[Theorem 3.5]{Mautner(p-adic.Lie)}),
      hence to any analytic group over any nonarchimedian field of
      characteristic $0$. 

\noindent(3)\quad
Any closed subgroup of a linear group over a local field with
respect
to \textbf{inner} automorphisms.
By the introductory remark it is enough to prove this for the group
$SL_n$, since we may embed any closed linear group therein as a
closed subgroup. 
We have already seen in Example \ref{1stEx:U,P}(1)
that a
contraction subgroup for an element of $SL_n$ is an algebraic       
subgroup, hence is closed. 
\end{bemerkung}\endrem

Further examples with all contraction groups for inner automorphisms
closed  may be obtained from the above list by forming projective limits
and restricted products.

\begin{itembemerkung}\label{G^+}
(1)\quad 
Note that item 1(c) above shows that all \tdlc\ MAP-groups
have
trivial inner contraction groups hence scale function identically
one (they are \emph{uniscalar}) if they are metric.
This is an improvement over the main result of \cite{lcMAP=unimod}.

\noindent(2)\quad 
Let $G^+$ be the closed subgroup generated by all
contraction groups of inner automorphisms.
It is stable under any (bicontinuous) automorphism of $G$,
hence is in particular normal.
As a consequence of Proposition~\ref{quotient-U}, 
for metric groups $G$, $G/G^+$ 
is uniscalar.
The scale function should thus characterize best the groups
satisfying $G=G^+$.
Every non-uniscalar topologically simple group
will belong to this class by Proposition~\ref{U bounded iff}, but
no solvable group will because
$G^+\subseteq \ol{[G,G]}$. 
This last remark suggests that we
define subgroups $G^{n+}$ by iterating the definition of $G^+$.
\end{itembemerkung}\endrem


\section{The tree-representation Theorem}\label{S:tree}
\label{tree}

We will show that, if $V$ is tidy for $\alpha$, then $V_{--}$ extended by
the cyclic group  generated by $\alpha$, that is, $V_{--} \rtimes \langle
\alpha \rangle$, has a representation onto some closed subgroup of the
automorphism group of a homogeneous tree.
This representation can be used to analyse groups of this type,
which are the simplest non-uniscalar groups.

Although this is not the approach we shall take, this representation is an
instance of a general construction in the Bass-Serre theory of group
actions on graphs which is described in
\cite{trees} and \cite{graphs}. 
Since $\alpha \colon V_- \to \alpha(V_-)$ is an injection into $V_-$, the group $V_{--} \rtimes \langle\alpha \rangle$ is isomorphic to the HNN extension
$$
V_- {\displaystyle
\mathop *_{\alpha(V_-)}} t, \text{ where }
t = \alpha^{-1} \colon \alpha(V_-)
\to V_-.
$$
\parbox{7.5cm}{The HNN extension is
defined in \cite[Example 3.5(v)]{graphs} to be the fundamental group of
the graph of groups at right. 
It is shown in \cite{trees} and \cite{graphs} how to represent HNN extensions on trees. 
For completeness, we shall describe the tree 
and the action of $V_{--}\rtimes \langle\alpha \rangle$ directly. 
Elements of $V_{--} \rtimes \langle \alpha
\rangle$ will be denoted $v\alpha^m$.}\quad
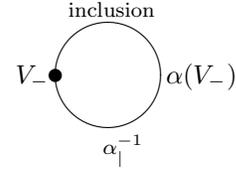
\begin{picture}(40,55)(0,20)
\put(35,30){\circle{50}}
\put(15,30){\circle*{5}}
\put(0,27){$V_-$}
\put(57,27){$\alpha(V_-)$}
\put(20,52){\scriptsize inclusion}
\put(33,0){\scriptsize $\alpha^{-1}_{|}$}
\end{picture}\\

The tree will be denoted by $T$. 
Its vertices are the left $V_-$-cosets in 
$V_{--} \rtimes \langle \alpha \rangle$. 
Distinct vertices $xV_-$ and $yV_-$ are linked by a directed edge 
$xV_-\to yV_-$ if and only if
\begin{equation}
\label{eq:edgedefn}
yV_- \subset xV_-\alpha = (x\alpha)\alpha^{-1}(V_-).
\end{equation}
Equivalently, there is an edge from $v\alpha^mV_-$ to $w\alpha^nV_-$ if
and only if $n = m+1$ and $w\in v\alpha^m(V_-)$. 
Since $\alpha^{-1}(V_-)$ is the union of $s(\alpha^{-1})$ $V_-$-cosets,
there  are $s(\alpha^{-1})$ out-edges from the vertex $xV_-$. 
Since $xV_- \subset \left((x\alpha^{-1})\alpha\right)\alpha^{-1}(V_-)$ 
there is one edge into $xV_-$, from the vertex $(x\alpha^{-1})V_-$. 
Hence each vertex in $T$ has degree $s(\alpha^{-1}) + 1$.

We show next that $T$ is a tree when $\alpha$ has infinite order. 
For each $n\in\mathbb{Z}$, denote the vertex $\alpha^nV_-$ by $V^{(n)}$. 
Then 
$V^{(n)} = \alpha^{n-1}\alpha(V_-)\alpha\subset \alpha^{n-1}V_-\alpha$ 
and it follows from (\ref{eq:edgedefn}) that
$$
\dots V^{(-2)}, V^{(-1)}, V^{(0)},
V^{(1)}, V^{(2)}, \dots
$$
is an infinite path, call it $P$, in $T$.
Consider a general vertex $v\alpha^mV_-$ in $T$, where $v\in
\alpha^{-n}(V_-)$. 
Then $v\alpha^mV_- = \alpha^mwV_-$, where 
$w =\alpha^{-m}(v) \in\alpha^{-(m+n)}(V_-)$. 
If $m+n\leq0$, it follows that $v\alpha^mV_- =V^{(m)}$. 
If $m+n > 0$, then $v\alpha^mV_-$ is descended by a path of
length $m+n$ from $V^{(-n)}$. 
In either case we see that $v\alpha^mV_-$ is connected to $P$. 
Therefore $T$ is connected. 
To see that there are no circuits in $T$, observe first of all that, 
since each vertex has only
one in-edge, any circuit 
without backtracking 
must be a directed path. 
However, the power of
$\alpha$ in each coset $xV_-$ strictly increases in the direction of any
path in $T$ and so there can be no directed circuit unless $\alpha$ has
finite order.

The action of $V_{--} \rtimes \langle \alpha \rangle$ on the vertices of
$T$ is the usual action of the group on a space of left $V_-$-cosets. 
It is clear that this action preserves the adjacency relation defined in
(\ref{eq:edgedefn}). 
Denote this action by 
$\rho \colon V_{--} \rtimes \langle\alpha \rangle \to \text{Aut}(T)$. 
Then $\rho$ is vertex and edge transitive on $T$ and so 
the quotient graph is a loop.

Denote the set of
ends of $T$ by $\partial T$; the end of the path $P$ corresponding to
$\{V^{(n)}\}_{n=1}^\infty$ by $\infty$; and the other end of $P$ by
$-\infty$. 
Then $\alpha$ acts as a translation of distance 1 on $P$ 
with $-\infty$ as the repelling end. 
Since every path in $P$ descends ultimately from $-\infty$
and it is the unique end of $T$ having this property, $-\infty$ is fixed
by $\rho(V_{--} \rtimes \langle\alpha \rangle)$. 
This may also be shown directly by verifying that
$v\alpha^mV^{(-r)} = V^{(m-r)}$ provided that $r$ is sufficiently large
that $\alpha^{r-m}(v)\in V_-$.

Recalling that the automorphism group of $T$ is itself a \tdlcG\ when
equipped with the topology of uniform convergence on compact sets, we are
in a position to state the tree-representation theorem.

\begin{theorem}
\label{Vtree}
Let $G$ be a \tdlcG, $\alpha$ an  automorphism of $G$
of infinite order and let $V$ be tidy for $\alpha$. 
Then $V_{--} \rtimes\langle\alpha \rangle \to \text{Aut}(T)$ is a continuous representation $\rho$ onto 
a closed subgroup of the automorphism group of a homogeneous tree
$T$ of degree $s(\alpha^{-1}) + 1$.
\begin{itemize}
\item[(1)] The action of $V_{--} \rtimes \langle \alpha \rangle$ on
$\text{Aut}(T)$: fixes an end, $-\infty$; is transitive on $\partial
T\setminus\{-\infty\}$; and the quotient graph is a loop.
\item[(2)] The stabiliser of each end in $\partial T\setminus\{-\infty\}$
is a conjugate of $V_0\rtimes \langle\alpha \rangle$. 
The kernel of $\rho$ is the largest compact normal 
$\alpha$-stable subgroup of $V_{--}$.
\item[(3)] The image of $V_{--}$ under $\rho$ is the set of elliptic
elements in $\rho(V_{--} \rtimes \langle \alpha \rangle)$.
\end{itemize}
\end{theorem}

\proof{}
Stabilisers of vertices $xV_-$ in $T$ form a subbasis for the topology
on $\text{Aut}(T)$. 
Since $\rho^{-1}(\text{stab}(xV_-)) = xV_-x^{-1}$,
which is open in $V_{--} \rtimes \langle \alpha \rangle$, it follows that
$\rho$ is continuous. 
The continuity of $\rho$ and the compactness of
$V_-$ imply that $\rho(V_{--} \rtimes \langle \alpha \rangle) \cap
\text{stab}(V^{(0)}) = \rho(V_-)$ is compact and hence closed. 
Therefore $\rho(V_{--} \rtimes \langle \alpha \rangle)$ is closed by \cite[II.(5.9)]{GMW115}. 

{``(1)'':}\quad
It remains only to show that $\rho$ is transitive on $\partial
T\setminus\{-\infty\}$ and, for this, it suffices to show that for
each $\omega\in \partial T\setminus\{-\infty\}$ there is $v\in
V_{--}$ such that $v\cdot\infty = \omega$. 
Let $\{v_n\alpha^nV_-\}_{n=1}^\infty$ be a path converging to $\omega$.
Since $(v_n\alpha^nV_-,v_{n+1}\alpha^{n+1}V_-)$ is an edge, it follows
from (\ref{eq:edgedefn}) that $v_{n+1}\alpha^{n+1}(V_-) \subset
v_n\alpha^n(V_-)$ for each $n$. 
Each of these sets is compact and so
$\bigcap_{n=1}^\infty v_n\alpha^n(V_-) \ne \emptyset$. 
Choose $v$ in this intersection, so that 
$v\alpha^n(V_-) = v_n\alpha^n(V_-)$ for each $n$. 
Then
$$
vV^{(n)} = v\alpha^n V_- = v\alpha^n(V_-)\alpha^n =
v_n\alpha^n(V_-)\alpha^n = v_n\alpha^nV_-
$$
for each $n$ and it follows that $v\cdot\infty = \omega$.

\noindent{``(2)'':}\quad
For the first part, it suffices to show that the stabiliser of $\infty$ equals $V_0\rtimes \langle \alpha \rangle$. 
It is clear that this group is contained in the stabiliser. 
Conversely, if $x$ stabilises $\infty$, then it leaves the path $P$ invariant and in fact translates $P$ by a distance, $d$ say. 
Then $x\alpha^{-d}$ fixes every vertex on $P$.
Hence $x\alpha^{-d} \in \bigcap_{n\in\ZZ} \alpha^n(V_-) = V_0$ and we have
$x\in V_0\rtimes \langle \alpha \rangle$.

The kernel of $\rho$ is the intersection of all the vertex stabilisers and
is therefore a compact normal $\alpha$-stable subgroup of $V_{--}$.
If $M$ is any such subgroup, then it fixes a point, $p$ say, in $T$
because it is compact. 
Since $M$ is also normal in $V_{--}$ and $\alpha$-invariant, 
it fixes every point in the $V_{--} \rtimes\langle \alpha \rangle$-orbit 
of $p$. 
Since $\rho$ is edge-transitive, it follows that $M$ fixes every point 
in $T$ and so $M$ is contained in the kernel of $\rho$.

\noindent``(3)'':
Since the image of $\rho$ fixes $-\infty$, if $g\in \im (\rho)$
fixes a point $p$, then it fixes every vertex on the path from $-\infty$
to $p$. 
Hence $g$ fixes a point if and only if there is  a
$k\in\mathbb{Z}$ such that $g.V^{(k)}=V^{(k)}$. 
The set of elliptics therefore coincides with
$$
\bigcup_{k\in\mathbb{Z}} \text{stab}(V^{(k)})\cap \im (\rho) =
\bigcup_{k\in\mathbb{Z}}
\rho\left(\alpha^k(V_-)\right) = \rho(V_{--}).
$$  
\qed

The representation $\rho$ restricts to give a representation of 
$\ol{U}_\alpha\rtimes \langle \alpha \rangle$ on $T$. 
It follows from Theorem~\ref{V--=} that, 
for $V_{--}$ metric, 
this is the same representation as obtained if
Theorem~\ref{Vtree} is applied with $G$ equal to $\ol{U}_\alpha$. 
Hence all the assertions of Theorem~\ref{Vtree} hold 
with $V_{--}$ replaced by $\ol{U}_\alpha$. 
However, more can be said about the representation of
$\ol{U}_\alpha\rtimes \langle \alpha \rangle$.

\begin{theorem}\label{T.alpha}
Let $G$ be a totally disconnected locally compact metric group 
and $\alpha$ an  automorphism of $G$ of infinite order. 
Let $\rho$ be the representation of 
$\ol{U}_\alpha\rtimes \langle \alpha \rangle$ 
on the tree $T$ as defined above.
\begin{itemize}
\item[(1)] The action of $U_\alpha$ is transitive on $\partial T
\setminus \{-\infty\}$ and is simply transitive if and only if $U_\alpha$
is closed.

\item[(2)] If $U_\alpha$ is closed, then $\rho$ is a topological
      isomorphism onto
      its image.
\end{itemize}
\end{theorem}
\proof{} 
{``(1)'':}\quad 
Theorem \ref{Vtree}(1) shows that $\ol{U}_\alpha$ is
transtitive on 
the set 
$\partial T \setminus \{-\infty\}$. 
Since, by Corollary
\ref{cor:Ualphabarfactor}, $\ol{U}_\alpha = U_\alpha U_0$ and since $U_0$
is the stabiliser of $\infty$, we have
$$
\partial T \setminus \{-\infty\} = \ol{U}_\alpha\cdot\infty =
U_\alpha\cdot\infty\,.
$$

The action is simply transitive if and only if $\text{stab}(\infty) =
U_0\cap U_\alpha$ is trivial. 
Since $U_0\cap U_\alpha$ is dense in $U_0$,
Theorem \ref{U0=1 iff} shows that $U_0\cap U_\alpha$ is trivial if and
only if $U_\alpha$ is closed.

\noindent{``(2)'':}\quad
If $U_\alpha$ is closed, then $\alpha$ is compactly contractive.
Hence  the compact kernel of $\rho$ is trivial and
$\rho$ is faithful. \qed

Subsection \ref{subsec:difference} examined the difference between the
groups $V_{--}$ and $U_\alpha$ where $\alpha$ is respectively shrinking
and contracting. 
Both groups act on $T$ with the $U_\alpha$-action being
the restriction of the $V_{--}$-action. The next example shows that the
difference between the groups may, or may not, be seen in the action on
the tree.

\begin{itembeispiel}
(1)\quad 
Let $G$ be the group of upper triangular matrices in 
$SL_2(\mathbb{Q}_p)$. 
The contraction group with respect to the automorphism 
$\alpha$ given by inner conjugation by 
$\text{diag}(p,p^{-1})$ 
is the group of strict upper triangular matrices. 
The subgroup of upper triangular matrices with entries in 
$\mathbb{Z}_p$ is compact open in $G$ and satisfies 
$\alpha(V)\subseteq V$ hence is tidy by Example~\ref{V=V-}. 
The group $V_{--}$ is then the subgroup of elements in $G$ 
having diagonal entries in $\mathbb{Z}_p^*$. 
The tree representation $\rho$ of $V_{--}$ is seen to be 
the composite of the inclusion in $SL_2(\mathbb{Q}_p)$ 
with the representation of $SL_2(\mathbb{Q}_p)$ on its 
Bruhat-Tits tree. 
Since $V_{--}$ properly contains $U_\alpha$, 
$\rho$ distinguishes between these groups. \\[-1ex]

\noindent(2)\quad
On the other hand, if we are given a closed contraction 
group $U_\alpha$ we may construct a strictly larger 
group $G$,  
an extension of $\alpha$ to $G$
and a tidy subgroup $V$ for $\alpha$ acting on $G$ 
with $G=V_{--}$ 
such that the tree representation 
$\rho$ defined by $V_-$ satisfies $\rho(V_{--})=\rho(U_\alpha)$ 
as follows. 
Take any nontrivial compact group $K$ and define $G$ by 
$G:=K\times U_\alpha$. 
Let $\alpha$ act trivially on $K$. 
Choose a tidy subgroup $O$ for $\alpha$ in $U_\alpha$. 
Then $V:=K\times O$ is tidy for $\alpha$ in $G$ and 
$V_{--}=K\times O_{--}=K\times U_\alpha=G$. 
Obviously $\rho(V_{--})=\rho(U_\alpha)$. 
\end{itembeispiel}

Part (2) of the example has $V_{--}$ as the direct product of $U_\alpha$
and the kernel of $\rho$. 
That is not special to this example. 
The kernel of $\rho$ and $\ol{U}_\alpha$ are both closed normal subgroups of $V_{--}$ and so their product is always a subgroup and it is closed because $\ker\rho$ is compact. 
Group isomorphism theorems plus the fact that, if
$U_\alpha$ is closed, then $U_\alpha\cap \ker\rho$ is trivial imply
the following result.

\begin{satz}
\label{prop:V--}
Let $V$ be tidy for the   
automorphism $\alpha$ of
the totally disconnected locally compact metric group $G$.\hfill\break
(1) If $\rho(\ol{U}_\alpha) = \rho(V_{--})$, then
$V_{--} = \ol{U}_\alpha \ker\rho$.\hfill\break
(2) If, furthermore, $U_\alpha$ is closed, then $V_{--} = {U}_\alpha
\times\ker\rho$.
\qed
\end{satz}

\begin{bemerkung}
Essentially the same tree as constructed here has been
associated with tidy subgroups and contraction groups elsewhere. 
In Theorem 3.4 of \cite{struc(tdlcG-graphs+permutations)} 
a tree construction is used to translate
the definition of tidy subgroup into permutation group theoretic terms.
The rooted tree constructed there is a branch of the tree constructed
above. 
Also, the underlying tree constructed in Theorem \ref{Vtree}, but
not the group representation $\rho$, is already implicit in
Proposition~3.7 in~\cite{contrAut(lcG)}, which studies the case where the
contraction group is closed.
\end{bemerkung}\endrem


\end{document}